\magnification \magstep1 \openup 2\jot
\def \qed {\vrule height6pt width6pt depth0pt}

\vbox{\vskip 1truecm}
\centerline{\bf{On the volume of the intersection of two $L_p^n$
balls }}

\centerline { by}

\centerline {G.~Schechtman\footnote { *}{ Supported in part by the US-
Israel BSF and by the Glikson Foundation} and
J.~Zinn\footnote {**}{ Supported in part by NSF DMS-86-01250
and by Texas Advanced Research Program}{\ }\footnote{\ }{
Grant no. 3825} }

{\  }

   {\  }

      {\  }

\beginsection 1. Introduction

This note deals with the following problem, the case $p=1$, $q=2$ of
which was introduced to us by Vitali Milman: What is the volume left
in the $L_p^n$ ball after removing a t-multiple of the $L_q^n$ ball?
Recall that the $L_r^n$ ball is the set $\{(t_1,t_2,\dots,t_n);\
t_i\in{\bf R},\ n^{-1}\sum_{i=1}^n|t_i|^r\le 1\}$ and note that for
$0<p<q<\infty$ the $L_q^n$ ball is contained in the $L_p^n$ ball.

In Corollary 4 below we show that, after normalizing Lebesgue measure
so that the volume of the $L_p^n$ ball is one, the answer to the
problem above is of order $e^{-ct^pn^{p/q}}$ for $T<t<{1\over 2}n^
{{1\over p}-{1\over q}}$, where $c$ and $T$ depend on $p$ and $q$ but
not on $n$.

The main theorem, Theorem 3, deals with the corresponding question
for the surface measure of the $L_p^n$ sphere. Theorem 3 and Corollary
4 together with some other remarks form Section 3. In Section 2
we introduce a class of random variables to be used in the proof of
the main theorem. These random variables are related to $L_p$ in the
same way that Gaussian variables are related to $L_2$.

\beginsection 2. Preliminaries

Here we introduce a class of random variables to be used in the proof of
the main theorem and summarize some of their properties.
Fix a $0<p<\infty$ and let\ $x,x_1,x_2,\dots,x_n$\
be independent random variables each with
density function \ $c_pe^{-t^p}$, $t>0$. Note that necessarily $c_p=
p/\Gamma (1/p)$. The first claim is known, though we could not locate a
reference.

{\ }

\noindent {\bf Lemma 1.} {\sl Put\
$S=\bigl( \sum_{i=1}^nx_i^p\bigr)^{1/p}$,\
then\  $\Bigl({{x_1}\over S},{{x_2}\over S},\dots,{{x_n}\over S}\Bigr)$
\ is uniformly distributed over the positive quadrant of the sphere
of \ $l_p^n$, i.e., over the
set\hfil\break
$\Delta_p = \{(t_1,t_2,\dots,t_n)\ ;\ t_i\ge~0,\ \sum t_i^p = 1\}$
equipped with the $(n-1)$-dimensional normalized Lebesgue measure.
Moreover,
$\Bigl({{x_1}\over S},{{x_2}\over S},\dots,{{x_n}\over S}\Bigr)$
is independent of $S$.}

{\ }

\noindent {\bf Proof.} For any Borel subset $A$ of $\Delta_p$,
$$\eqalign{
P\Bigl(\Bigl(&{{x_1}\over S},{{x_2}\over S},\dots,{{x_n}\over S}\Bigr)
\in A \Big|S=a\Bigr)=\cr
&= \lim_{\epsilon\rightarrow 0}{{P((x_1,\dots,x_n)\in
{\bf R}_+A\ \&\ a-\epsilon\le S\le a+\epsilon )}\over {P(a-\epsilon\le
S \le a+\epsilon)}}\cr
&=\lim_{\epsilon\rightarrow 0}
\int_{(t_1,\dots,t_n)\in{\bf R}_+A\atop
(a-\epsilon)^p<\sum t_i^p<(a+\epsilon
)^p}\  e^{-\sum t_i^p}dt\Big/
\int_{(t_1,\dots,t_n)\in{\bf R}_+^n\atop
(a-\epsilon)^p<\sum t_i^p<(a+\epsilon)^p}\
e^{-\sum t_i^p}dt\cr
&\le\limsup_{\epsilon\rightarrow 0}e^{-(a-\epsilon)^p+(a-\epsilon)^p}
\Bigr)^p
\int_{(t_1,\dots,t_n)\in{\bf R}_+A\atop
(a-\epsilon)^p<\sum t_i^p<(a+\epsilon
)^p}\ dt\Big/
\int_{(t_1,\dots,t_n)\in{\bf R}_+^n\atop
(a-\epsilon)^p<\sum t_i^p<(a+\epsilon)^p}\ dt\cr
&=\lambda(A),\cr }$$
where $\lambda$ is the normalized Lebesgue
measure on $\Delta$. Similarly,
$$P\Bigl(\Bigl({{x_1}\over S},{{x_2}\over S},\dots,{{x_n}\over S}\Bigr)
\in A \Big|S=a\Bigr)\ge\lambda(A).$$
This proves that
$P\Bigl(\Bigl({{x_1}\over S},{{x_2}\over S},\dots,{{x_n}\over S}\Bigr)
\in A\Bigr)
=\lambda(A)$ and that
$\Bigl({{x_1}\over S},{{x_2}\over S},\dots,{{x_n}\over S}\Bigr)$
is independent of $S$.\hfill\qed

{\ }

In the next claim we gather some more properties of the random variables
$x_i$.

{\ }

\noindent{\bf Lemma 2.} {\sl Let $x,x_1,\dots,x_n$ be as above, then
\item {1.} $c_p$ is bounded away from zero and infinity when
$p\rightarrow\infty$.
\item {2.} For all $h>0$ and all $0<p<\infty$,
${\bf E}e^{-hx^p}=\Bigl({1\over {1+h}}\Bigr)^{1/p}$.
In particular,$${\bf E}e^{-hx^p}\ge e^{-h/p}\ \ for \ all\ \  h>0\ \  and
\ \ {\bf E}e^{-hx^p}\le e^{-h/2p}\ \  for \ all\ \  0<h\le1.$$
\item {3.} For all $0<u<\infty$ and all $0<p<\infty$, $P(x^p>u)\ge
{{c_p}\over {2p}}e^{-2u}$. If $p\ge1$ and $u\ge 1$,
 then also $P(x^p>u)\le
{{c_p}\over p}e^{-u/2}$. In particular, for $p\ge 1$ and all $u$,
$P(x^p>u)\le Ce^{-u/2}$ for some universal $C$.
\item {4.} For all $1\le p\le q<\infty$,
 ${\bf E}
\Bigl(\sum_{i=1}^nx_i^q\Bigr)^{1/q}$ is equivalent, with universal
constants, to\  $q^{1/p}n^{1/q}$, if $q\le\log n$, and to
$(\log n)^{1/p}$ otherwise.}

{\ }

\noindent{\bf Proof.} 1. Follows easily from the fact that
$c_p=p/\Gamma (1/p)=\Gamma({1\over p}+1)^{-1}$.

2. is a simple computation.

3. is also simple, here is a sketch of the proof.
$$\eqalign{P(x^p>u)&=c_p\int_{u^{1/p}}^\infty e^{-t^p}dt\cr
&\ge
c_p\int_{u^{1/p}}^{(u+1)^{1/p}} {{pt^{(p-1)}}\over {p(u+1)^{(p-1)/p}}}
\,e^{-t^p}dt\cr
&={{c_p}\over {p(u+1)^{(p-1)/p}}}\Bigl(1-{1\over e}\Bigr) e^{-u}\cr
&\ge{{c_p}\over {2p(u+1)}}\, e^{-u}\cr&\ge
{{c_p}\over {2p}}\,e^{-2u}.\cr}$$
The other inequality in 3 is proved in a similar way.

4. First note that for all $0<p,q<\infty$
$${\bf E}x^q=c_p\int_0^\infty t^qe^{-t^p}dt={{c_p}\over p}\Gamma
\Bigl({{q+1}\over p}\Bigr)$$
so that, by the triangle inequality and 1, if $1\le p\le q <\infty$
$${\bf E}\Bigl(\sum_{i=1}^nx_i^q\Bigr)^{1/q}\le
\Bigl(\sum_{i=1}^n{\bf E}x_i^q\Bigr)^{1/q}=
\Bigl({{c_p}\over p}\Gamma
\Bigl({{q+1}\over p}\Bigr)\Bigr)^{1/q}n^{1/q}\le Cq^{1/p}n^{1/q}$$
for some universal $C$.
For the lower bound in the case $q\le\log n$,
divide $\{1,2,\dots,n\}$ into approximately $n/e^q$
disjoint sets of
cardinality approximately $e^q$ each, then
$$\eqalign{{\bf E}\Bigl(\sum_{i=1}^nx_i^q\Bigr)^{1/q}&=
{\bf E}\Bigl(\sum_j\Bigl(\sum_{i\in\sigma_j}x_i^q\Bigr)^{q/q}\Bigr)^{1/q}
\cr
&\ge
{\bf E}\Bigl(\sum_j\bigl(\max_{i\in\sigma_j}x_i\bigr)^q
\Bigr)^{1/q}\cr
&\ge
\Bigl(\sum_j\bigl({\bf E}\max_{i\in\sigma_j}x_i\bigr)^q
\Bigr)^{1/q}\cr
&\ge c^\prime(\log e^q)^{1/p}(n/e^q)^{1/q}\cr
&\ge c^{\prime\prime}q^{1/p}n^{1/q}.\cr}$$

Now, for the case $q>\log n$ we note first that, by 3,
$$P(\max_{1\le i\le n}x_i>t)\ge 1-\bigl(1-{{c_p}\over
{2p}}e^{-2t^p}\bigr)^n.$$
For $n$ smaller than an absolute multiple of $p$, the lower bound
follows easily from the fact that ${\bf E}x_1$ is larger that a universal
positive constant, so assume that $n\ge 20p/c_p$ and put
$t=2^{-1/p}\Bigl(\log{{nc_p}\over {2p}}\Bigr)^{1/p}$.Then, for some
universal $c$,
$$P(\max_{1\le i\le n}x_i > c(\log n)^{1/p})\ge 1/2.$$
In particular, ${\bf E}\max_{1\le i\le n}x_i\ge c(\log n)^{1/p}$,
which implies the lower bound in this case since
$\Bigl(\sum_{i=1}^nx_i^q\Bigr)^{1/q}$ is universally equivalent to
$\max_{1\le i\le n}x_i .$
The upper bound in this case, though a bit harder, is also standard
and since we don't use it in the sequel we shall leave it to the reader.
\hfill\qed

 {\ }

The statement in 4, for the case $p=2$, was noticed by the first named
author several years ago while seeking a precise estimate for the
dimension of the Euclidean sections of $l_p^n$ spaces (see
\lbrack MS\rbrack\ p.145 Remark 5.7).
The original proof was more complicated.
The proof presented here is an adaptation of a proof of the case $p=2$
shown to us by J. Bourgain.

\beginsection 3. The main result

\noindent{\bf Theorem 3.} {\sl
For all $1\le p<q<\infty$ there are constants
$c=c(p,q)$ , $C=C(p,q)$ and $T=T(p,q)$ such that if $\mu$ denotes the
normalized Lebesgue measure on the positive quadrant of the unit sphere
of $L_p^n$ then
$$\mu(\Vert u\Vert_{L_q^n} > t) \le \exp(-ct^pn^{p/q})\eqno{(1)}$$
for all $t > T$, and
$$\mu(\Vert u\Vert_{L_q^n} > t) \ge \exp(-Ct^pn^{p/q})\eqno{(2)}$$
for all $2\le t\le {1\over 2}n^{{1\over p}-{1\over q}}$.

\noindent Moreover, for
$q>2p$ (or any other universal positive multiple of $p$), one can take
$c(p,q)={\gamma\over p}$,
 $C(p,q)={\Gamma\over p}$  and
$T(p,q)=\tau \min\{ q,log n\}^{1/p}\le q^{1/p}$.
Here $\gamma$ ,$\Gamma$ and $\tau$ are universal constants.}

{\ }

\noindent{\bf Proof.} By Lemma 1 above,
$$\mu(\Vert u\Vert_{L_q^n} > t) = P(n^{{1\over p}-{1\over q}}
{(\sum_{i=1}^n x_i^q)^{1/q}}/
{(\sum_{i=1}^n x_i^p)^{1/p}}>t)$$
where $x_i$ are independent random variables each with density
$c_p e^{-t^p}$.
Assume, for the simplicity of the presentation, that $n$ is even.
Put\  $S=(\sum_{i=1}^n x_i^p)^{1/p}$ \
and let $p_j$,\ \ \  $j=1,2,\dots,n/2$\  be positive
numbers with sum $\le 1/2$. Then
$$\eqalignno{ P\bigl(n^{{1\over p}-{1\over q}}
&{(\sum_{i=1}^n x_i^q)^{1/q}}/
{(\sum_{i=1}^n x_i^p)^{1/p}}>t\bigr) =\cr
&= P\Bigl(\sum_{i=1}^n x_i^q >
{{t^q(\sum_{i=1}^n x_i^p)^{q/p}}\over {n^{{q\over p }-1}} }\Bigr)\cr
&\le\sum_{i=1}^{n/2}P(x_i^* >
tp_i^{1/q}S/n^{{1\over p} -{1\over q}})
+ P\bigl(\sum_{i={n\over 2} + 1}^n x_i^{*q} >
t^qS^q/2n^{{q\over p} -1} \Bigr)&{(3)}\cr}$$
where $\{ x_j^*\}$ denotes the nonincreasing rearrangement of
$\{|x_j|\}$.

Since $$\eqalign{\sum_{j={n\over 2}+1}^n x_i^{*q}&\le
{n\over 2}\,x_{n\over 2}^{*q} \le
{n\over 2} ({2\over n}\sum_{i=1}^{n/2}x_i^{*p})^{q/p}\cr
&\le 2^{{q\over p}-1}S^q/n^{{q\over p}-1} \ \ ,}$$
we get that, if $t\ge 2^{1/p}$, the second term in (3) is zero.

To evaluate the first term in (3), fix $1\le j\le n/2$. Then,
$$\eqalign{P(x_j^* > tp_j^{1/q}S/n^{{1\over p}-{1\over q}})
&\le {n\choose j}P(x_1,\dots,x_j
 > tp_j^{1/q}S/n^{{1\over p}-{1\over q}})
\cr &\le {n\choose j}P\bigl(x_1^p,\dots,x_j^p >
t^pp_j^{p/q}\sum_{i=j+1}^n x_i^p/n^{1-{p/q}}\bigr).}$$

From Lemma 2 (first 3 and then 2) we get that the last
expression is dominated by
$${n\choose j}C^j{\bf E}
\exp\bigl(-jp_j^{p/q}t^p\sum_{i=j+1}^n x_i^p /n^{1-{p/q}}\bigr)
 $$
$$\le{n\choose j}C^j
\exp(-jp_j^{p/q}t^p(n-j)/2pn^{1-{p/q}})$$
for some universal $C$. Note that
the last inequality holds if $jn^{{p/q}-1}p_j^{p/q}t^p \le 1$.
If this is not the case the probability we are trying to evaluate
is zero. Finally, the last term is dominated by
$$\exp\Bigl(j\bigl(log{en\over j} + C-
{{p_j^{p/q}t^pn^{p/q}}\over {4p}}\bigr)\Bigr).\eqno{(4)}$$
Now, for $\alpha$ to be chosen momentarily, let $p_j$, $j=1,\dots,n/2$,
be such that
$$j\bigl(log{en\over j} + C-
{{p_j^{p/q}t^pn^{p/q}}\over {4p}}\bigr)= - \alpha n^{p/q}t^p$$
i.e.,
$$p_j = \bigl( 4p{{log{en\over j}}\over {t^pn^{p/q}}} +
{{4Cp}\over {t^pn^{p/q}}} + \alpha {{4p\over j}}\bigr)^{q/p}.$$
We thus get that,for some universal constant $C$,
$$p_j \le 2^{{q\over p}-1} {{(Cp)^{q/p}(log {{en}\over j})^{q/p}}\over
{t^qn}} + 2^{{q\over p}-1}
\alpha^{q/p}{{(4p)^{q/p}}\over {j^{q/p}}}.\eqno{(5)}$$
It is easy to see that, for $1\le p<q<\infty$,
$$\sum_{j=1}^{n/2}(log{{en}\over j})^{q/p}\le
An\min\{q^{q/p},(log n)^{q/p}\}$$
for some universal $A$.
Thus the sum over $j$ of the first
terms in (5) is smaller than 1/4 if, for some universal
$\gamma$,\  $t>
\gamma \min\{q^{1/p},(\log n)^{1/p}\}$.
The sum over $j$ of the second terms in (5) is bounded by 1/4 if
$\alpha < B{1\over p}({q\over p}-1)^{p/q}$, for some universal $B$.
Choosing $\alpha$ to satisfy this inequality and using (3),(4) and (5)
we get that, for $t > \gamma \min\{q^{1/p},(logn)^{1/p}\}$,
$$\mu(\Vert u\Vert_{L_q^n} > t) \le {n\over 2} e^{-\alpha n^{p/q}t^p}.$$
Under the conditions on $t$, the factor $n/2$ can be absorbed in the
second term (changing $\alpha$ to another constant of the same order of
magnitude as a function of $p$),
 thus proving (1).

We now turn to the proof of the lower bound (2) which is simpler.
Using Claim 1 again,
$$\eqalign{\mu\Vert u\Vert_{L_q^n} > t) &= P(n^{{1\over p}-{1\over q}}
{(\sum_{i=1}^n x_i^q)^{1/q}}/
{(\sum_{i=1}^n x_i^p)^{1/p}}>t)\cr
&\ge P(x_1 > St/n^{{1\over p}-{1\over q}})\cr
&=P\Bigl(x_1>
{t\over {(n^{(1-p/q)}-t^p)^{1/p}}}
\bigl(\sum_{i=2}^nx_i^p\bigr)^{1/p}\Bigr).\cr}$$
Since $t^p\le {1\over 2}n^{(1-p/q)}$, this dominates
$$P\Bigl(x_1>
{{2^{1/p}t}\over {n^{{1\over p}-{1\over q}}}}
\bigl(\sum_{i=2}^nx_i^p\bigr)^{1/p}\Bigr).$$
Now, by Claim 2.3.,
$$\eqalignno{
P\Bigl(x_1>
{{2^{1/p}t}\over {n^{{1\over p}-{1\over q}}}}
\bigl(\sum_{i=2}^nx_i^p\bigr)^{1/p}\Bigr)
&\ge
{{c_p}\over {2p}}{\bf E}\exp(-4t^p\sum_{i=2}^nx_i^p/n^{(1-p/q)})\cr
&={{c_p}\over {2p}}\bigl(
{\bf E}\exp(-4t^px_1^p/n^{(1-p/q)})\bigr)^{n-1}\cr
&={{c_p}\over {2p}}\Bigl({1\over{1+{{4t^p}\over {n^{1-p/q}}}}}
\Bigr)^{(n-1)/p}&{(by\  Claim\ 2.2.)}\cr
&\ge {{c_p}\over{2p}}\exp\Bigl(-{{4t^p(n-1)}\over{pn^{(1-p/q)}}}\Bigr)\cr
&\ge{{c_p}\over{2p}}e^{  4t^pn^{p/q}/p}.}$$
Finally observe that, since $c_p$ is bounded away from zero and $t\ge 2$,
the factor ${{c_p}\over {2p}}$ can be absorbed in the second term
(changing $4$ to another universal constant).\hfill\qed

{\ }

\noindent{\bf Remarks:} 

1. It follows from the proof that, for $n$ large enough and $q$
close to $p$, one can take $c(p,q)={c\over p}\Bigl({q\over p}-1\Bigr)$
for some universal constant $c$.

2. It follows from the statement of the theorem that, for $q=\infty$,
$$\mu(\Vert u\Vert_\infty > t) \le e^{-\gamma t^p/p}$$
for all $t > \tau(\log n)^{1/p}$, and
$$\mu(\Vert u\Vert_\infty > t) \ge e^{-\Gamma t^p/p}$$
for all $2\le t\le {1\over 2}n^{1\over p}$, where $\gamma$, $\Gamma$ and
$\tau$ are universal constants.

3. Note that it follows from Claim 1 and Claim 2.4. that the order of
magnitude of $T$ is the correct one.

4. The restriction $p\ge 1$ in Theorem 3 above and in Corollary 4 below
can be replaced by $p>0$ if one replaces the inequality $t\ge 2$ with
$t\ge d$, for some $d$ depending only on $p$ and $q$, and removes the
``moreover" part. We didn't check the dependence of the constants on
$p$ and $q$ in this case.

The last remark is that one can get a similar statement for the full
balls. We state it as a corollary.

{\ }

\noindent{\bf Corollary 4.} {\sl
For all $1\le p<q<\infty$ there are constants
$c=c(p,q)$ , $C=C(p,q)$ and $T=T(p,q)$ such that if $\nu$ denotes the
normalized Lebesgue measure on the ball
of $L_p^n$ then, for all $n$ large enough,
$$\nu(\Vert u\Vert_{L_q^n} > t) \le \exp(-ct^pn^{p/q})\eqno{(6)}$$
for all $t > T$, and
$$\nu(\Vert u\Vert_{L_q^n} > t) \ge \exp(-Ct^pn^{p/q})\eqno{(7)}$$
for all $2\le t\le {1\over 2}n^{{1\over p}-{1\over q}}$.
\noindent Moreover, for
$q>2p$ (or any other universal positive multiple of $p$), one can take
$c(p,q)={\gamma\over p}$ ,
 $C(p,q)={\Gamma\over p}$  and
$T(p,q)=\tau \min\{ q,log n\}^{1/p}\le q^{1/p}$,
where $\gamma$ ,$\Gamma$ and $\tau$ are universal constants.}

  {\ }

The proof follows easily from Theorem 3 and the formula
$$\nu(A)=n\int_0^1r^{n-1}\mu\Bigl({A\over r}\Bigr)dr$$
which holds for all Borel sets $A$ in the ball of $L_p^n$.

{\ }

\beginsection References

\item{\lbrack MS\rbrack}V. D. Milman and G. Schechtman,
Asymptotic Theory  of Finite Dimensional Normed Spaces,
~~Lecture Notes in Math. Vol. 1200, Springer
(1986).

{\baselineskip=14pt

{\  }

{\ }
 {\obeylines \parindent=0pt \lineskip=2pt

Department of Theoretical Mathematics,
The Weizmann Institute of Science,
Rehovot, Israel

{\ }

Department of Mathematics,
Texas A\&M University,
College Station, Texas 77843, USA
{\ }}}
\vfill
\end